\documentclass[12pt]{article}
\usepackage{latexsym}
\usepackage{amsfonts}
\usepackage{amssymb}
\usepackage{amscd}
\usepackage{array}
\usepackage{amsmath}
\usepackage{epsfig}
\begin{document}
\title{\textbf{On the Hard Lefschetz property of stringy Hodge numbers}\\ \date{}}
\author{Jan Schepers\footnote{Supported by the Netherlands Organisation for
Scientific Research (NWO).}}
\maketitle
\begin{center}
\footnotesize{\textbf{Abstract}}
\end{center}
{\footnotesize For projective varieties with a certain class of `mild' isolated
singularities and for projective threefolds with arbitrary Gorenstein canonical
singularities, we show that the stringy Hodge numbers satisfy the Hard Lefschetz
property (i.e.\ $h_{st}^{p,q}\leq h_{st}^{p+1,q+1}$ for $p+q \leq d-2$, where $d$ is
the dimension of the variety). This result fits nicely with a 6-dimensional
counterexample of Musta\c{t}\u{a} and Payne for the Hard Lefschetz property for
stringy Hodge numbers in general. We also give such an example, ours is a
hypersurface singularity.}

\section{Introduction} 
\noindent \textbf{1.1.} Stringy Hodge numbers of projective varieties with
Gorenstein canonical singularities were introduced by Batyrev in \cite{Batyrev}.
They are defined if the \textit{stringy $E$-function}, which is in general a
rational function of two variables $u$ and $v$, is in fact a polynomial. The idea is
that they should be the Hodge numbers of a conjectural `string cohomology'. Several
constructions of such string cohomology spaces were made by Borisov and Mavlyutov in
\cite{BorisovMavlyutov}, and they also made a connection to the orbifold cohomology
of Chen and Ruan from \cite{ChenRuan}. Moreover, Yasuda showed that the stringy
Hodge numbers are the Hodge numbers of the orbifold cohomology for varieties with
Gorenstein quotient singularities (see \cite[Remark 1.4 (2)]{Yasuda}). \\

\noindent \textbf{1.2.} In this paper we want to study under which conditions the
Hard Lefschetz property holds for stringy Hodge numbers. By the Hard Lefschetz
property we mean the inequalities that are imposed if there would be an analogue of
the Hard Lefschetz Theorem for the conjectural string cohomology, i.e.
\[ h_{st}^{p,q} \leq h_{st}^{p+1,q+1} \text{ for } p+q \leq d-2,    \] where $d$ is
the dimension of the variety.
Fernandez gave in \cite{Fernandez} a criterion for the Hard Lefschetz Theorem to
hold for orbifold cohomology. He also shows that this criterion fails in general for
generic Calabi-Yau hypersurfaces in the weighted projective space
$\mathbb{P}(1,1,1,3,3)$ (\cite[Example 4.4]{Fernandez}). However, it follows for
instance from the main theorem below that this cannot be seen from the Hodge or
Betti numbers of the orbifold cohomology ($h_{orb}^{p,q}=h_{st}^{p,q}\leq
h_{orb}^{p+1,q+1}=h_{st}^{p+1,q+1}$ for $p+q \leq 1$). In fact, the orbifold Hodge
numbers where in this case already studied by Poddar (see \cite{Poddar}, Section 5
and especially Corollary 2). \\

\noindent \textbf{1.3.} Musta\c{t}\u{a} and Payne first gave an example of a
6-dimensional projective toric variety with stringy Betti number $b_{6,st} =
h_{st}^{3,3}$ strictly smaller than $b_{4,st}=h^{2,2}_{st}$ (see \cite[Example
1.1]{MustataPayne}). This example was used to disprove a conjecture of Hibi on the
unimodality of the so-called $h$-vector or $\delta$-vector of a reflexive polytope.
This vector actually gives the stringy Betti numbers of the toric variety defined by
the fan over the faces of the polytope (\cite[Theorem 7.2]{BatyrevDais} and
\cite[Theorem 1.2]{MustataPayne}). For many more examples of reflexive polytopes
with non-unimodal $h$-vector we refer to \cite{Payne}.\\

\noindent \textbf{1.4.} The main result of this paper is the following. The used
notions are explained in Section 2.\\

\noindent \textbf{Main theorem.} \textsl{Let $Y$ be either
\begin{itemize}
\item a projective variety of dimension $d=3$ with Gorenstein canonical
singularities, or
\item a projective variety of dimension $d\geq 4$ with at most isolated Gorenstein
singularities that admit a log resolution with all discrepancy coefficients of
exceptional components $> \lfloor \frac{d-4}{2} \rfloor$. (This condition does not
depend on the chosen log resolution.)
\end{itemize}
Write the stringy $E$-function of $Y$ as a power series $\sum_{i,j\geq 0}
b_{i,j}u^iv^j$. Then for $ i+j\leq d-2$, we have $(-1)^{i+j}b_{i,j}\leq
(-1)^{i+j+2}b_{i+1,j+1}$. In particular, if the stringy $E$-function of $Y$ is a
polynomial, then $h_{st}^{p,q}(Y)\leq h_{st}^{p+1,q+1}(Y)$ for $p+q\leq d-2$.}\\

\noindent For the proof of this theorem we refer to Section 3. In Section 4 we
compare the above theorem with the example of Musta\c{t}\u{a} and Payne. We also
discuss an explicit example of a 6-dimensional projective variety with an isolated
canonical hypersurface singularity that does not satisfy the Hard Lefschetz
property.\\

\noindent \textbf{Acknowledgements.} I want to thank Wim Veys for helpful
discussions and for encouraging me to communicate these results. I also want to
thank the referee of an earlier paper for bringing references \cite{Fernandez} and
\cite{MustataPayne} to my attention.

\section{Stringy Hodge numbers} 
\noindent \textbf{2.1.} Let $X$ be an arbitrary complex algebraic set of dimension
$d$. It is well known that the cohomology with compact support $H_c^{\bullet}(X)$
carries a natural mixed Hodge structure (we always use cohomology with coefficients
in the complex numbers). The data of this mixed Hodge structure are put in the
Hodge-Deligne polynomial
\[ H(X;u,v):= \sum_{p,q=0}^d \left[\sum_{i=0}^{2d} (-1)^i h^{p,q}(H_c^i(X))\right]
u^pv^q,  \] where $h^{p,q}$ denotes the dimension of the $(p,q)$-component
$H^{p,q}(H_c^i(X))$. The Hodge-Deligne polynomial is a generalized Euler
characteristic: if $Y$ is a Zariski-closed subset of $X$, then
$H(X;u,v)=H(X\setminus Y;u,v) + H(Y;u,v)$ and for a product of algebraic sets $X$
and $X'$ we have $H(X\times X';u,v) = H(X;u,v) H(X';u,v)$. \\

\noindent \textbf{2.2.} Now let $Y$ be a normal irreducible variety. Assume that $Y$
is $\mathbb{Q}$-Gorenstein; i.e.\ $rK_Y$ is a Cartier divisor for some $r\in
\mathbb{Z}_{>0}$. If $K_Y$ is already Cartier then $Y$ is called Gorenstein. For
example, all hypersurfaces and more generally all complete intersections are
Gorenstein. Choose a log resolution $f:X\to Y$ of $Y$. This is a proper birational
map $f$ from a nonsingular variety $X$ such that the exceptional locus is a divisor
$D$ with only smooth components $D_i,i\in I,$ that have normal crossings. We can
write 
\[ rK_X - f^*(rK_Y) = \sum_{i\in I} b_iD_i, \] with all $b_i\in \mathbb{Z}$. Using
$\mathbb{Q}$-coefficients this becomes $K_X-f^*(K_Y) = \sum a_i D_i$ with $a_i =
b_i/r$. The rational number $a_i$ is called the discrepancy coefficient of $D_i$. We
call $Y$ log terminal, canonical or terminal if all $a_i > -1, \geq 0$ or $> 0$
respectively (this does not depend on the chosen log resolution). \\

\noindent \textbf{2.3. Definition} (\cite[Definition 3.1]{Batyrev})\textbf{.} Let
$Y$ be log terminal. Choose a log resolution $f:X\to Y$ with irreducible exceptional
components $D_i,i\in I$. Denote the discrepancy coefficient of $D_i$ by $a_i$. For a
subset $J$ of $I$ we use the notations
\[ D_J:= \bigcap_{j\in J} D_j \quad \text{ and }\quad D_J^{\circ}:= D_J \setminus
\bigcup_{i\in I\setminus J} D_i. \] This gives a stratification of $X$ as
$\coprod_{J\subset I} D_J^{\circ}$. The \textit{stringy $E$-function} of $Y$ is
defined by the formula
\[ E_{st}(Y;u,v) := \sum_{J\subset I} H(D_J^{\circ};u,v) \prod_{j\in J}
\frac{uv-1}{(uv)^{a_j+1}-1}.\] 
Batyrev used motivic integration to show that this formula does not depend on the
chosen log resolution (\cite[Theorem 3.4]{Batyrev}).\\

\noindent \textbf{2.4. Remark.}
\begin{enumerate}
\item If $Y$ is Gorenstein (and thus automatically canonical) then $E_{st}(Y;u,v)$
is a rational function in $u$ and $v$. It lives in $\mathbb{Z}[[u,v]] \cap
\mathbb{Q}(u,v)$.
\item If $Y$ is smooth, then $E_{st}(Y;u,v)=H(Y;u,v)$. If $Y$ admits a crepant
resolution (i.e.\ a log resolution $f:X\to Y$ such that $K_X = f^*(K_Y)$) then
$E_{st}(Y;u,v) = H(X;u,v)$. More generally, for a projective birational morphism
$g:Y' \to Y$ from a normal variety $Y'$ such that $K_Y' = g^*(K_Y)$ one has
$E_{st}(Y';u,v)=E_{st}(Y;u,v)$ (\cite[Theorem 3.12]{Batyrev}).
\item An alternative formula for $E_{st}(Y;u,v)$ is
\[ E_{st}(Y;u,v) = \sum_{J\subset I} H(D_J;u,v) \prod_{j\in J}
\frac{uv-(uv)^{a_j+1}}{(uv)^{a_j+1}-1}. \]
\end{enumerate}

\noindent \textbf{2.5.} Assume now that $Y$ is in addition projective of dimension
$d$. Batyrev proves the following relation in \cite[Theorem 3.7]{Batyrev}:
\[ E_{st}(Y;u,v) = (uv)^d E_{st}(Y;u^{-1},v^{-1}). \tag{1} \] If $Y$ is also
Gorenstein canonical and if $E_{st}(Y;u,v)$ is a polynomial $\sum_{p,q}
a_{p,q}u^pv^q$ then Batyrev defines the \textit{stringy Hodge numbers} of $Y$ as
$h_{st}^{p,q}(Y):= (-1)^{p+q} a_{p,q}$. Remark the following:
\begin{enumerate}
\item If $Y$ is smooth then the stringy Hodge numbers are equal to the usual Hodge
numbers of $Y$. If $Y$ has a crepant desingularization $X$, then the stringy Hodge
numbers of $Y$ are equal to the Hodge numbers of $X$.
\item $h_{st}^{p,q}(Y) = h_{st}^{q,p}(Y)$.
\item From (1) above it follows that $h_{st}^{p,q}(Y)$ can only be nonzero for $0
\leq p,q \leq d$ and that $h_{st}^{p,q}(Y)=h_{st}^{d-p,d-q}(Y)$. 
\item From Remark 2.4 (3) we have $h_{st}^{0,0}(Y)=1$. 
\end{enumerate}
Batyrev also made the following very intriguing conjecture.\\

\noindent \textbf{Conjecture} (\cite[Conjecture 3.10]{Batyrev})\textbf{.} Stringy
Hodge numbers are nonnegative.\\

\noindent \textbf{2.6. Example.} Canonical surface singularities are classified and
are precisely the $A$-$D$-$E$ singularities. It is well known that they admit a
crepant resolution and thus the conjecture above is trivially true for them by
Remark 2.5 (1). Of course, the Hard Lefschetz property is also satisfied in this
case.\\

\noindent In \cite{SchepersVeys} the conjecture was proved for the same class of
varieties that is treated by the main theorem of this paper, and thus also for
threefolds in full generality.

\section{Proof of the main theorem}
\noindent \textbf{3.1.} Let us for convenience repeat the statement of the theorem.\\

\noindent \textbf{Theorem.} \textsl{Let $Y$ be either
\begin{itemize}
\item a projective variety of dimension $d=3$ with Gorenstein canonical
singularities, or
\item a projective variety of dimension $d\geq 4$ with at most isolated Gorenstein
singularities that admit a log resolution with all discrepancy coefficients of
exceptional components $> \lfloor \frac{d-4}{2} \rfloor$.
\end{itemize}
Write the stringy $E$-function of $Y$ as a power series $\sum_{i,j\geq 0}
b_{i,j}u^iv^j$. Then for $ i+j\leq d-2$, we have $(-1)^{i+j}b_{i,j}\leq
(-1)^{i+j+2}b_{i+1,j+1}$.}\\ 

\noindent \textbf{Remark.}
\begin{enumerate}
\item In particular, if the stringy $E$-function of $Y$ is a polynomial, then
$h_{st}^{p,q}(Y)\leq h_{st}^{p+1,q+1}(Y)$ for $p+q\leq d-2$.
\item It is not hard to check that the lower bound on the discrepancies for the
second class of varieties from the theorem does not depend on the chosen log
resolution. Note that isolated terminal four- and fivefold singularities are
included in the theorem.
\item Note also that $b_{0,0}=E_{st}(Y;0,0)=1$ by Remark 2.4 (3), so $b_{i,i} \geq
1$ for $ i \leq d/2$.
\end{enumerate}

\vspace{0.3cm}

\noindent \textit{Proof of the theorem.} Let us first treat the second case. So $Y$
is of dimension $d\geq 4$. Let $f:X \to Y$ be a log resolution with $X$ projective
and such that $f$ is an isomorphism when restricted to the inverse image of the
nonsingular part of $Y$. Denote by $D$ the total exceptional locus of all singular
points. In \cite[Remark 3.5 (3)]{SchepersVeys} the following description of the
numbers $(-1)^{i+j}b_{i,j}$ for $i+j\leq d$ was given in this case:
\[  (-1)^{i+j}b_{i,j} = \dim \ker \bigl(H^{d-i,d-j}(H^{2d-i-j}(X))\to
H^{d-i,d-j}(H^{2d-i-j}(D))\bigr) + S_{i,j},       \]
where one has to remark that
\begin{enumerate}
\item the map $H^k(X)\to H^k(D)$ induced by inclusion is surjective for $k\geq d$
and thus $H^k(D)$ carries then a pure Hodge structure,
\item $S_{i,j}$ is nonnegative and for even $d$ only nonzero for $i=j=d/2$ and for
odd $d$ only nonzero for $i=j=(d-1)/2$ and $\{i,j\}=\{(d-1)/2,(d+1)/2\}$. The term
$S_{i,j}$ has to be introduced for contributions of the exceptional components with
the lowest allowed discrepancy coefficients. If one does not put the above lower
bound on the discrepancy coefficients, it becomes more often nonzero, much more
complicated and can even be negative (see the example in \cite{SchepersVeys2}).
\end{enumerate}
If we denote $\ker \bigl(H^{d-i,d-j}(H^{2d-i-j}(X))\to
H^{d-i,d-j}(H^{2d-i-j}(D))\bigr)$ by $K^{2d-i-j}_{i,j}$ then it suffices to prove
that $\dim K^{2d-i-j}_{i,j} \leq \dim K^{2d-i-j-2}_{i+1,j+1}$ for $i+j \leq d-2$.
Denote $\ker (H^k(X)\to H^k(D))$ by $K^k$.\\

\noindent We will use the following construction of de Cataldo and Migliorini. Embed
$Y$ in a projective space $\mathbb{P}^r$ and take a generic hyperplane section $Y_s$
of $Y$. So $Y_s$ is nonsingular and does not contain any of the singular points of
$Y$. Let $X_s:=f^{-1}(Y_s)$ and denote by $\eta_s$ the fundamental class of $X_s$ in
$H^{1,1}(H^2(X))$. If one dualizes the surjective map $H^k(X) \to H^k(D)$ for $k\geq
d$ and uses Poincar\'e duality on $X$, one obtains an injection $H_k(D) \to
H^{2d-k}(X)$. Now the spaces
\[ H^0(X),H^1(X), \frac{H^2(X)}{H_{2d-2}(D)}, \ldots, 
\frac{H^{d-1}(X)}{H_{d+1}(D)}, H^d(X),\]
\[K^{d+1},\ldots,K^{2d-2},H^{2d-1}(X),H^{2d}(X) \] satisfy the Hard Lefschetz
Theorem with respect to the cup product with $\eta_s$. This result is discussed by
de Cataldo and Migliorini in Section 2.3 and the beginning of Section 2.4 from
\cite{dCM2} for dimensions 3 and 4, but it is not hard to see that their argument
works in any dimension. It also follows from their earlier work \cite{dCM1}, see
Section 2.4 there. We note that the choice of the space $H^d(X)$ in the middle is
somewhat arbitrary, it can be replaced by any subspace containing the image of
$\frac{H^{d-2}(X)}{H_{d+2}(D)}$. We can take $K^d$ for that. To prove this, it
suffices to show that \[ H^{d-2}(X)\begin{array}{c} \cup \eta_s \\ \longrightarrow
\\ \ \end{array} H^d(X) \to H^d(D) \vspace{-0,5cm} \] forms a complex. If we
dualize, this means that \[ H_d(D)\to H_d(X)\to H_{d-2}(X) \] should be a complex as
well, 
 where $H_d(X)\to H_{d-2}(X)$ corresponds to intersecting with $X_s$. And this is
clear.\\

\noindent Summarizing, we obtain that the maps $\cup \eta_s: K^{d+k} \to K^{d+k+2}$
are surjective for $0\leq k \leq d-2$. Since $\cup \eta_s$ is a morphism of Hodge
structures of type (1,1), we also get that $\cup \eta_s: K^{2d-i-j-2}_{i+1,j+1} \to
K^{2d-i-j}_{i,j}$ is surjective for $i+j\leq d-2$.\\

\noindent Now let $Y$ be a projective threefold with arbitrary Gorenstein
singularities. By the main theorem of \cite{Reid2} we can find a projective variety
$Z$ with terminal singularities and a projective birational crepant morphism $g:Z\to
Y$. So $E_{st}(Z)=E_{st}(Y)$ by Remark 2.4 (2). The point is that terminal threefold
singularities are automatically isolated (see for instance \cite[Corollary
4-6-6]{Matsuki}) and thus we can apply the above reasoning for $Z$ (the used results
of de Cataldo and Migliorini remain valid for $Z$, as well as the description of the
numbers $(-1)^{i+j}b_{i,j}$ from \cite{SchepersVeys} for $i+j\leq 3$, now with
$S_{i,j}$ always zero). \hfill $\blacksquare$

\section{Examples}
\noindent \textbf{4.1. Example.} We first compare Example 1.1 of \cite{MustataPayne} with the main theorem. Let $f$ be the vector
$(\frac{1}{3},\frac{1}{3},\frac{1}{3},\frac{1}{3},\frac{1}{3},\frac{1}{3})$ in
$\mathbb{R}^6$ and $N$ the lattice $\mathbb{Z}^6 +  \mathbb{Z}\cdot f$. Denote the
standard basis vectors of $\mathbb{R}^6$ by $e_1,\ldots,e_6$. Musta\c{t}\u{a} and
Payne consider the polytope $P$ with vertices
\[  \{e_1,\ldots , e_6 , e_1-f , \ldots, e_6-f\}. \] It is a reflexive polytope and
the projective toric variety $Y$ defined by the fan $\Sigma$ over the faces of the
polytope has stringy $E$-function
\[ (uv)^6 + 6(uv)^5 + 8(uv)^4 + 6(uv)^3 + 8(uv)^2 + 6uv +1.\]
The fan $\Sigma$ has eight cones of maximal dimension, namely 
\begin{itemize}
\item a cone $\sigma$ generated by $e_1,\ldots , e_6$,
\item $\tau$ generated by $e_1-f, \ldots , e_6-f$,
\item and six non-simplicial cones $\rho_i$ generated by all vectors to vertices
except $e_i$ and $e_i-f$.
\end{itemize}
As sketched in \cite[Example 3.1]{MustataPayne} one can make a triangulation of the boundary of $P$ by taking the convex hulls of $\{e_1, \ldots, \widehat{e_j}, \ldots, e_k , e_k-f, \ldots,e_6-f\}$ and $\{e_1, \ldots , e_j ,e_j-f , \ldots , \widehat{e_k - f } , \ldots , e_6\}$ for $1\leq j < k \leq 6$. This triangulation is regular in the sense of \cite[Section I.1.F]{BrunsGubeladze} as can be seen for instance by repeatedly applying Lemma 1.72 from that book. This implies that the toric variety $Z$ given by the fan over this triangulation is projective. Moreover, the toric morphism from $Z$ to $Y$ is crepant, since no new rays and hence no exceptional components of codimension 1 where introduced. Thus $E_{st}(Z) = E_{st}(Y)$. Since the singular locus of $Z$ is given by the union of the orbits of the torus action corresponding to cones in the fan that cannot be generated by a part of the basis of the lattice, we see that $Z$ has exactly two isolated singular points coming from the cones $\sigma$ and $\tau$. The resolution of singularities of these points is particularly easy: we subdivide the fan by adding the rays generated by $f$ and by $-f$. This does introduce two irreducible exceptional components of codimension 1 and from the theory of toric varieties it is well known that there discrepancy coefficients are 1. So the example of the variety $Z$ shows that the lower bound on the discrepancies in the main theorem is crucial. Even for isolated singularities with a very easy resolution, the theorem cannot be extended. \\

\noindent \textbf{4.2.} We want to conclude this paper by discussing another
6-dimensional example that was obtained independently of the one of Musta\c{t}\u{a}
and Payne. The resolution of singularities is much more complicated, but it has the
advantage of being a hypersurface singularity. We will need the formula for the
Hodge-Deligne polynomial of a Fermat hypersurface and for a quasihomogeneous affine hypersurface with an isolated singularity at the origin.\\

\noindent Denote the Fermat hypersurface of dimension $d$ and degree $l$ by $Y_l^{(d)}$. So $Y_l^{(d)}$ is given by
\[ \{x_0^l + \cdots + x_{d+1}^l =0\} \subset \mathbb{P}^{d+1}_{\mathbb{C}}.    \]
Dais shows in \cite[Lemma
3.3]{Dais} that the Hodge-Deligne polynomial of $Y^{(d)}_l$ is given by 
\[ H(Y^{(d)}_l;u,v) := \sum_{p=0}^{d} u^p\left(v^p +
(-1)^d\mathcal{G}(d+1,p+1\,|\,l-1,p)v^{d-p} \right), \]
where \[ \mathcal{G}(\kappa,\lambda\,|\,\nu,\xi) :=
\sum_{j=0}^{\lambda} (-1)^j {\kappa + 1 \choose j} {\nu(\lambda -j) + \xi \choose
\kappa} \] for $(\kappa,\lambda,\nu,\xi)\in \mathbb{Z}_{\geq 0}^4$ and $\kappa \geq
\lambda$ (if $m>n$, the binomial coefficient ${n\choose m}$ must be interpreted as
0).\\

\noindent Let $f\in \mathbb{C}[x_1,\ldots,x_{r+1}]$ be a
quasihomogeneous polynomial of degree $d$ with respect to the weights 
$w_1,\ldots,w_{r+1}$ and assume that the origin is an isolated singularity of
$Y:=f^{-1}(0)$. According to
\cite[Section 2]{Dais} the Hodge-Deligne polynomial of $Y$ equals
\[(uv)^r +(-1)^{r-1}(uv-1)\sum_{p=0}^{r-1}
\dim_{\mathbb{C}}M(f)_{(p+1)d-(w_1+\cdots+w_{r+1})}u^pv^{r-1-p}, \]
where $M(f)_{(p+1)d-(w_1+\cdots+w_{r+1})}$ denotes the piece of degree $(p+1)d-(w_1+\cdots+w_{r+1})$ of the Milnor algebra 
\[ M(f):=\frac{\mathbb{C}[x_1,\ldots,x_{r+1}]}{\left(\frac{\partial f}{\partial
x_1},\ldots,\frac{\partial f}{\partial x_{r+1}} \right)}.      \] 
Indeed, this is a graded $\mathbb{C}$-algebra if we give $x_i$
degree $w_i$. The needed dimensions can be computed from the Poincar\'e series 
\[ P_{M(f)}(t)
:= \sum_{k\geq 0} (\dim_{\mathbb{C}}M(f)_k)t^k,     \]
which in this case simply equals 
\[ P_{M(f)}(t)= \frac{(1-t^{d-w_1})\cdots (1-t^{d-w_{r+1}})}{(1-t^{w_1})\cdots
(1-t^{w_{r+1}})}.  \]

\vspace{0,5cm}

\noindent \textbf{4.3. Example.} It is not so easy to give the variety $Y$ whose stringy $E$-function we want to
compute. Let us start from 
\[ Y' := \{ x_1^5z^3 + x_2^5z^3 + x_3^8 +x_4^8+x_5^8+x_6^8+x_7^8 = 0\} \subset
\mathbb{P}^7,\]
where we consider $z=0$ as the hyperplane at infinity. At infinity there is a
singular $\mathbb{P}^1$ and the origin of the affine chart $z\neq 0$ is
singular as well. Our variety $Y$ will consist of a resolution of the singular
$\mathbb{P}^1$ and will thus have one isolated hypersurface singularity. Let us first
describe the resolution process at infinity. Thereby we want to compute the
Hodge-Deligne polynomial of the nonsingular part $Y_{ns}$ of $Y$, since the stringy $E$-function of $Y$ can be written as $H(Y_{ns};u,v)\,+$\,contribution of the singular point. So we keep track
of the contributions in every step. We blow up in the singular line. The exceptional
locus after this first step consists of five disjoint components $D_1^{\infty},
\ldots, D_5^{\infty}$ (all isomorphic to $\mathbb{P}^5$) and one other
component $D_6^{\infty}$, also isomorphic to $\mathbb{P}^5$ and
singular for the strict transform of $Y'$. The intersection of $D_6^{\infty}$ and
another $D_i^{\infty}$ is isomorphic to $\mathbb{P}^4$. Since we are only
interested in the nonsingular part and since we will blow up in $D_6^{\infty}$ in the
following step, the contribution of the first step to the Hodge-Deligne polynomial
of $Y_{ns}$ will be $5(uv)^5$ (the nonsingular exceptional part are five disjoint
$\mathbb{A}^5$'s). So in the second step we blow up with
$D_6^{\infty}$ as center. Again there appear five new disjoint components
$E_1^{\infty},\ldots,E_5^{\infty}$ and one component $E_6^{\infty}$ intersecting
the others transversally and singular for the strict transform of $Y'$. Here the
$E_i^{\infty}$ are $\mathbb{P}^1$-bundles over
$\mathbb{P}^4$ and $E_6^{\infty}$ is isomorphic to
$\mathbb{P}^5$. The intersection of $E_6^{\infty}$ with an
$E_i^{\infty}$ is isomorphic to $\mathbb{P}^4$ and thus the
contribution of this step to the Hodge-Deligne polynomial of $Y_{ns}$ is
$5(uv)^5+5(uv)^4+5(uv)^3+5(uv)^2+5uv$. In the next step we blow up with
$E_6^{\infty}$ as center. There are two new exceptional components. The first, $F^{\infty}$, is a
$\mathbb{P}^1$-bundle over $\mathbb{P}^4$ that contains
five disjoint singular components for the strict transform of $Y'$, all isomorphic
to $\mathbb{P}^4$. The second component $G^{\infty}$ is isomorphic to a
$\mathbb{P}^1\times \mathbb{P}^1$-bundle over the Fermat
hypersurface $Y_8^{(3)}$ (notation as in 4.2). The intersection $F^{\infty}\cap
G^{\infty}$ is isomorphic to a $\mathbb{P}^1$-bundle over $Y_8^{(3)}$. All of this
means that the contribution of this step to the Hodge-Deligne polynomial of
$Y_{ns}$ equals
\[ (uv-4)((uv)^4+(uv)^3+(uv)^2+uv+1)  + ((uv)^2+uv)H(Y_8^{(3)};u,v).\]
In the final step we blow up in the five remaining singular components. This gives
five new disjoint exceptional components $H_1^{\infty},\ldots, H_5^{\infty}$, whose
Hodge-Deligne polynomial equal
\[ (uv)^5+2(uv)^4+2(uv)^3+2(uv)^2+2uv+1+uv H(Y_8^{(3)};u,v).\]
So the total contribution of the exceptional locus above the singular
$\mathbb{P}^1$ of $Y'$ to the Hodge-Deligne polynomial of $Y_{ns}$ is
\[16(uv)^5+12(uv)^4+12(uv)^3+12(uv)^2+12uv+1+((uv)^2+6uv)H(Y_8^{(3)};u,v).\]
To compute the Hodge-Deligne polynomial of $Y_{ns}$ we must add the Hodge-Deligne
polynomial of the nonsingular part of $Y'$. At infinity this is the double
projective cone over $Y_8^{(3)}$ minus the singular $\mathbb{P}^1$,
with contribution 
\[(uv)^2H(Y_8^{(3)};u,v),\]
and the contribution of the nonsingular part of $Y'$ in the affine chart $z\neq 0$
can be computed by the method of Dais from 4.2; it equals
\[ (uv)^6-1 - (uv-1)(140u^4v+140uv^4+4060u^2v^3+4060u^3v^2).\]
The formula for $H(Y_8^{(3)};u,v)$ is 
\[ (uv)^3+(uv)^2+uv+1-35u^3-35v^3-1015u^2v-1015uv^2,\]
so finally the Hodge-Deligne polynomial of $Y_{ns}$ equals
\begin{equation*} \begin{split} E :=&\ (uv)^6 +
18(uv)^5+20(uv)^4+20(uv)^3+20(uv)^2+18uv-210u^5v^2\\ &
-210u^2v^5-6090u^4v^3-6090u^3v^4
-70u^4v-70uv^4-2030u^3v^2\\ &-2030u^2v^3.
\end{split} \end{equation*}

\vspace{0,5cm}

\noindent Next we compute the contribution of the singular point given by the origin of
\[ \{ x_1^5+ x_2^5+x_3^8+x_4^8+x_5^8+x_6^8+x_7^8 =0 \}\subset
\mathbb{A}^7\]
to the stringy $E$-function of $Y$ by a log resolution. We remark that one can also use the combinatorial procedure of \cite[Section 4]{SchepersVeys2}. We first blow up in the
singular point itself. This gives five exceptional components $D_1,\ldots,D_5$, all
isomorphic to $\mathbb{P}^5$ and intersecting in the new singular locus
(isomorphic to $\mathbb{P}^4$). Blowing up in this intersection gives
two new exceptional components, but one of them (isomorphic to a
$\mathbb{P}^1$-bundle over $\mathbb{P}^4$) is singular for
the strict transform of $Y$. The other one, called $E$, is isomorphic to a
$\mathbb{P}^2$-bundle over $Y_8^{(3)}$. Its intersection with the $D_i$
is covered by the new singular locus. Blowing up in this singular locus creates six
new exceptional components. Five of them behave similarly. They are isomorphic to a
$\mathbb{P}^1$-bundle over $\mathbb{P}^4$, they are
disjoint and each of them has an intersection with one $D_i$ isomorphic to
$\mathbb{P}^4$ (they are now the only components that intersect the
$D_i$). We call these components $F_1,\ldots,F_5$ and choose the numbering
compatible with those of the $D_i$. The sixth exceptional component is called $G$
and is also isomorphic to a $\mathbb{P}^1$-bundle over
$\mathbb{P}^4$. The new singular locus is $E\cap G$ and it is
isomorphic to a $\mathbb{P}^1$-bundle over $Y_8^{(3)}$. The intersection of an
$F_i$ and $G$ is isomorphic to $\mathbb{P}^4$ and has a
$Y_8^{(3)}$ in common with the singular locus and the intersection of $E$ and an
$F_i$ is only 3-dimensional and is covered by the singular locus. In the final step
we blow up in the remaining singular locus. There is one new exceptional component,
called $H$. Its Hodge-Deligne polynomial equals
$((uv)^2+7uv + 1)H(Y_8^{(3)};u,v)$ and it splits off $E$ from the other components.
The intersections of $H$ with the $F_i$, with $G$ and with $E$ are all isomorphic to
a $\mathbb{P}^1$-bundle over $Y_8^{(3)}$. This final blow up also adds
$uvH(Y_8^{(3)};u,v)$ to the Hodge-Deligne polynomial of the $F_i$. And there are
threefold intersections $F_i\cap G\cap H$ isomorphic to $Y_8^{(3)}$. The result is a
normal crossing divisor with the following intersection diagram:
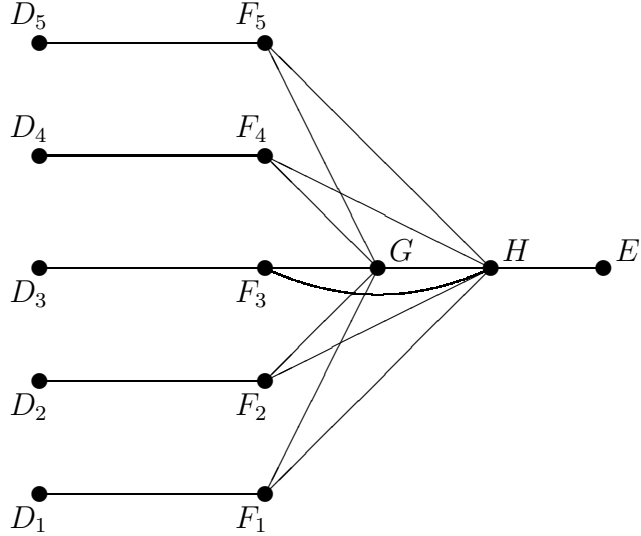
\begin{center}
\begin{picture}(102,70)
\put(20,5){\circle*{2}} \put(20,20){\circle*{2}} \put(20,35){\circle*{2}}
\put(20,50){\circle*{2}} \put(20,65){\circle*{2}} \put(50,35){\circle*{2}} 
\put(50,5){\circle*{2}} \put(50,20){\circle*{2}} \put(50,50){\circle*{2}}
\put(50,65){\circle*{2}} \put(65,35){\circle*{2}} \put(80,35){\circle*{2}}
\put(95,35){\circle*{2}} \put(20,35){\line(1,0){75}}  \put(20,65){\line(1,0){30}}
\put(20,50){\line(1,0){30}} \put(20,20){\line(1,0){30}} \put(20,5){\line(1,0){30}} 
\qbezier(50,35)(65,28)(80,35) \put(16,0.5){$D_1$} \put(16,15.5){$D_2$}
\put(16,30.5){$D_3$} \put(16,52.5){$D_4$} \put(16,67.5){$D_5$}  \put(66.5,36){$G$} 
\put(81.5,36){$H$} \put(96.5,36){$E$} \put(65,35){\line(-1,1){15}}
\put(65,35){\line(-1,-1){15}} \put(65,35){\line(-1,2){15}}
\put(65,35){\line(-1,-2){15}} \put(80,35){\line(-1,1){30}}
\put(80,35){\line(-1,-1){30}} \put(80,35){\line(-2,1){30}}
\put(80,35){\line(-2,-1){30}} \put(46,0.5){$F_1$} \put(46,15.5){$F_2$}
\put(46,30.5){$F_3$} \put(46,52.5){$F_4$} \put(46,67.5)
 {$F_5$}
\end{picture}
\end{center}
The discrepancy coefficients are 1 for the $D_i$ and 0 for all the other components.
Thus one can compute that
\begin{equation*} \begin{split}  F:=&\ 
8(uv)^5+15(uv)^4+10(uv)^3+15(uv)^2+8uv+1-70u^5v^2-70u^2v^5\\ &
-2030u^4v^3-2030u^3v^4
-210u^4v-210uv^4-6090u^3v^2-6090u^2v^3.
\end{split} \end{equation*} is the contribution of the singular point to the stringy $E$-function.
Then the stringy $E$-function of $Y$ is just $E+F$ which equals
\begin{equation*} \begin{split}   E_{st}(Y;u,v) =& \  
(uv)^6+26(uv)^5+35(uv)^4+30(uv)^3+35(uv)^2+26uv+1\\ &
-280u^5v^2-280u^2v^5-8120u^4v^3-8120u^3v^4
-280u^4v \\ &-280uv^4 -8120u^3v^2-8120u^2v^3,
\end{split} \end{equation*}
and thus $h_{st}^{2,2}(Y) > h_{st}^{3,3}(Y)$.\\

\footnotesize{

$\phantom{some place}$

\noindent Jan Schepers\\ Universiteit Leiden\\ Mathematisch Instituut\\ Niels
Bohrweg 1 \\ 2333 CA Leiden\\ The Netherlands\\ \emph{E-mail}:
jschepers@math.leidenuniv.nl

\end{document}